\documentclass[12pt]{article}
\usepackage{amssymb}
\usepackage{amsmath}

\begin{document}

\title{Monstrous and Generalized Moonshine and Permutation Orbifolds}
\author{Michael P. Tuite \\
School of Mathematics, Statistics and Applied Mathematics, \\
National University of Ireland, \\
Galway, Ireland.}
\maketitle

\begin{abstract}
We consider the application of permutation orbifold constructions towards a
new possible understanding of the genus zero property in Monstrous and
Generalized Moonshine. We describe a theory of twisted Hecke operators in
this setting and conjecture on the form of Generalized Moonshine replication
formulas.
\end{abstract}

\section{Introduction}

The Conway and Norton Monstrous Moonshine Conjectures \cite{CN}, the
construction of the Moonshine Module \cite{FLM} as an orbifold vertex
operator algebra and the completion of the proof of Monstrous Moonshine by
Borcherds \cite{Bo2} provided much of the motivation for the development of
Vertex Operator Algebras (VOAs) e.g. \cite{Bo1}, \cite{FLM}, \cite{Ka}, \cite%
{MN}. Another highlight of VOA theory is Zhu's study of the modular
properties of the partition function (and\ $n$-point functions) for generic
classes of VOAs \cite{Z}. Zhu's ideas were generalized to include orbifold
VOAs \cite{DLM} whose relevance to Monstrous Moonshine is emphasized in
refs. \cite{T1}, \cite{T2}. Norton's Generalized Moonshine Conjectures \cite%
{N2} concerning centralizers of the Monster group has yet to be generally
proven using either Borcherds' approach or orbifold partition function
methods although some progress has recently been made in refs. \cite{H} and 
\cite{T3}, \cite{IT1}, \cite{IT2} respectively.

In this note, we sketch a possible new approach to these areas based on
permutation orbifold VOA constructions \cite{DMVV}, \cite{BDM}. In
particular, we introduce a theory of twisted Hecke operators generalizing
classical Hecke operators in number theory\ e.g. \cite{Se}. We then discuss
permutation orbifold constructions where the classical Hecke operators
naturally appear and finite group and permutation orbifold constructions
where the twisted Hecke operators appear. Using these ideas we formulate a
conjecture on the nature of Generalized Moonshine replication formulas
generalizing replication formulas for the classical $J$ function and
McKay-Thompson series in Monstrous Moonshine. Detailed proofs will appear
elsewhere \cite{T4}.

\section{Replication Formula for the $J$ Function}

We begin with a brief review of Faber polynomials, Hecke algebras and the
replication formula for the classical $J$ function. Consider 
\begin{equation*}
t(q)=q^{-1}+0+\sum\limits_{k\geq 1}a(k)q^{k},
\end{equation*}%
the formal series in $q$. Define the Faber polynomial $P_{n}(x)$ for $t(q)$
to be the unique $n^{th}$ order polynomial with coefficients in $\mathbb{Z}%
[a(1),\ldots ,a(n-1)]$ such that 
\begin{equation}
P_{n}(t(q))=q^{-n}+O(q).  \label{Pnf}
\end{equation}%
Thus $P_{1}(x)=x$, $P_{2}(x)=x^{2}-2a(1)$, $P_{3}(x)=x^{3}-3a(1)x-3a(2)$
etc. The Faber polynomials for $t$ satisfy the following generating relation
e.g. \cite{Cu}, \cite{N1} 
\begin{equation}
\exp (-\sum\limits_{n\geq 1}\frac{p^{n}}{n}P_{n}(x))=p(t(p)-x),
\label{Faber gen}
\end{equation}%
for formal parameter $p$.

Let $f(\tau )$ be a meromorphic function of $\tau \in \mathbb{H}$, the upper
half complex plane. Then for integer $k\geq 1$ define a right modular group
action on $f$ for $\gamma =\left( 
\begin{array}{cc}
a & b \\ 
c & d%
\end{array}%
\right) \in \Gamma =SL(2,\mathbb{Z})$ as follows%
\begin{equation}
(f|_{k}\gamma )(\tau )=(c\tau +d)^{-k}f(\gamma \tau ),  \label{g_gamma}
\end{equation}%
with $\gamma \tau =\frac{a\tau +b}{c\tau +d}$. Then $f$ is a modular form of
(necessarily even) weight $k$ if $f$ is holomorphic in $\tau $ and 
\begin{equation*}
f|_{k}\gamma =f.
\end{equation*}%
Define the standard Hecke operators $T(n)$ for $n\geq 1$ with the following
action on a modular form $f$ of weight $k$ \cite{Se} 
\begin{equation}
T(n)f(\tau )=\frac{1}{n}\sum\limits_{a\geq 1,ad=n}a^{k}\sum\limits_{0\leq
b<d}f(\frac{a\tau +b}{d}).  \label{Tnf}
\end{equation}%
These satisfy the Hecke algebra 
\begin{eqnarray}
T(mn) &=&T(m)T(n),\text{ \ \ \ }(m,n)=1,  \notag \\
T(p)T(p^{m}) &=&T(p^{m+1})+p^{k-1}T(p^{m-1}),\text{ \ \ \ }m\geq 1\text{, \
prime }p\text{ .}  \label{Hecke_algebra}
\end{eqnarray}%
One finds%
\begin{equation}
(T(n)f)|_{k}\gamma =T(n)f,  \label{Tnf_gam}
\end{equation}%
i.e. $T(n)f$ is also a modular form of weight $k$.

The classical example is the Eisenstein series $G_{k}$ of even weight $k\geq
4$ (with $G_{k}(\tau )=0$ for odd $k$) [op.cit.] 
\begin{eqnarray}
G_{k}(\tau ) &=&\sum_{\substack{ m,n\in \mathbb{Z}  \\ (m,n)\neq (0,0)}}%
\frac{1}{(m\tau +n)^{k}}  \notag \\
&=&2\zeta (k)+2\frac{(2\pi i)^{k}}{(k-1)!}\sum_{n\geq 1}\sigma
_{k-1}(n)q^{n},  \label{Ek}
\end{eqnarray}%
where $q=\exp (2\pi i\tau )$ and $\sigma _{k}(n)=\sum_{d|n}d^{k}$.
Furthermore, $G_{k}$ is an eigenfunction of $T(n)$ with eigenvalue
determined by the coefficient of $q^{n}$ normalized to the coefficient of $q$
i.e. 
\begin{equation}
T(n)G_{k}=\sigma _{k-1}(n)G_{k}.  \label{Tn_Ek}
\end{equation}%
Thus it follows from (\ref{Hecke_algebra}) that for $k$ odd 
\begin{eqnarray}
\sigma _{k}(mn) &=&\sigma _{k}(m)\sigma _{k}(n),\text{ \ \ \ }(m,n)=1, 
\notag \\
\sigma _{k}(p)\sigma _{k}(p^{m}) &=&\sigma _{k}(p^{m+1})+p^{k}\sigma
_{k}(p^{m-1}),\text{ \ \ \ }m\geq 1\text{, \ prime }p\text{ .}
\label{sigma_fact}
\end{eqnarray}%
In fact, it is easy to check directly that (\ref{sigma_fact}) holds for all $%
k\in \mathbb{C}$.

The classical modular invariant function of weight $0$ is given by 
\begin{eqnarray*}
J(\tau ) &=&1728\frac{G_{4}^{3}}{G_{4}^{3}-G_{6}^{2}}-744 \\
&=&\sum\limits_{k\in \mathbb{Z}}c(k)q^{k}=q^{-1}+0+196884q+21493760q^{2}+%
\ldots
\end{eqnarray*}%
with standard normalization $c(-1)=1$ and $c(0)=0$. $J$ is a hauptmodul for
the genus zero group $\Gamma $ and is thus a generator for the field of
modular invariants. Thus 
\begin{eqnarray}
T(n)J(\tau ) &=&\sum\limits_{a\geq 1,a|n}\frac{1}{a}\sum\limits_{s\in 
\mathbb{Z}}c(\frac{ns}{a})q^{as}  \label{TnJ_sum} \\
&=&\frac{1}{n}q^{-n}+O(q).  \label{TnJ_lead}
\end{eqnarray}%
is a polynomial in $J$ which from (\ref{Pnf}) and (\ref{TnJ_lead}) must be 
\begin{equation}
T(n)J(\tau )=\frac{1}{n}P_{n}(J(\tau )),  \label{Hecke_J}
\end{equation}%
where $P_{n}$ is the Faber polynomial corresponding\ to $J$. Eqn. (\ref%
{Hecke_J}) is called the replication formula for $J$.

Eqn. (\ref{TnJ_sum}) also implies that 
\begin{equation*}
\sum\limits_{n\geq 1}p^{n}T(n)J(\tau )=\sum\limits_{r\geq 1,s\in \mathbb{Z}%
}c(rs)\sum\limits_{a\geq 1}\frac{1}{a}p^{ar}q^{as}=-\sum\limits_{r\geq
1,s\in \mathbb{Z}}c(rs)\log (1-p^{r}q^{s}).
\end{equation*}%
Then (\ref{Faber gen}) implies the famous $J$ function denominator formula 
\cite{N1}, \cite{Bo2} 
\begin{eqnarray}
\exp (-\sum\limits_{n\geq 1}p^{n}T(n)J(\tau )) &=&\prod\limits_{r\geq 1,s\in 
\mathbb{Z}}(1-p^{r}q^{s})^{c(rs)}  \notag \\
&=&p(J(p)-J(q)).  \label{denom_formula}
\end{eqnarray}%
This formula is one of the cornerstones of Borcherds' celebrated proof of
the genus zero Moonshine property where (\ref{denom_formula}) is a
denominator formula for a particular generalized Kac-Moody algebra
constructed from the Moonshine Module $V^{\natural }$ \cite{Bo2}.

\section{Twisted Hecke Operators and Eisenstein Series}

The definitions of modular functions and Hecke operators above can be
generalized to "twisted" versions as follows. We define a twisted modular
form of integer weight $k$ to be a holomorphic (in $\tau $) function $%
f=f((\theta ,\phi ),\tau )$ for $(\theta ,\phi )\in U(1)\times U(1)$ such
that%
\begin{equation*}
f|_{k}\gamma =f,
\end{equation*}%
where%
\begin{equation*}
(f|_{k}\gamma )((\theta ,\phi ),\tau )=(c\tau +d)^{-k}f(\gamma (\theta ,\phi
),\gamma \tau ),
\end{equation*}%
with left group action 
\begin{equation}
\gamma (\theta ,\phi )=(\theta ^{a}\phi ^{b},\theta ^{c}\phi ^{d}).
\label{gamma_the_phi}
\end{equation}%
Clearly the case $(\theta ,\phi )=(1,1)$ defines a standard modular form of
weight $k$.

We can extend the definition of the Hecke operator $T(n)$ to twisted modular
forms as follows: 
\begin{equation}
T(n)f((\theta ,\phi ),\tau )=\frac{1}{n}\sum\limits_{a\geq
1,ad=n}a^{k}\sum\limits_{0\leq b<d}\ f((\theta ^{a}\phi ^{b},\phi ^{d}),%
\frac{a\tau +b}{d}),  \label{Tnf_twist}
\end{equation}%
which includes the standard definition in the case $(\theta ,\phi )=(1,1)$.
For $\phi =1$ and $\theta ^{m}=1$ for integer $m$, this Hecke operator is
essentially that which appears in Borcherds' proof \cite{Bo2} and is
discussed at length in ref. \cite{F}.

We also define a homothety operator\footnote{%
A similar operator is defined in the standard case \cite{Se}.} 
\begin{equation}
R(n)f((\theta ,\phi ),\tau )=\ f((\theta ^{n},\phi ^{n}),\tau ).
\label{Rn_twist}
\end{equation}%
These operators satisfy the Hecke algebra \cite{T4} 
\begin{eqnarray}
R(mn) &=&R(m)R(n)  \notag \\
R(m)T(n) &=&T(n)R(m)  \notag \\
T(mn) &=&T(m)T(n),\text{ \ \ \ }(m,n)=1,  \notag \\
T(p)T(p^{m}) &=&T(p^{m+1})+p^{k-1}T(p^{m-1})R(p),\text{ \ \ \ }m\geq 1\text{%
, \ prime }p\text{,}  \label{Twisted_Hecke}
\end{eqnarray}%
and one again finds 
\begin{equation}
(T(n)f)|_{k}\gamma =T(n)f,  \label{Tnf_gam_twist}
\end{equation}%
i.e. $T(n)f$ is also a twisted modular form of weight $k$.

A twisted Eisenstein series $G_{k}((\theta ,\phi ),\tau ))$ of weight $k\geq
1$ can also be defined \cite{DLM}, \cite{MTZ}. In particular, for $k\geq 4$
we define\footnote{%
The notation used here differs from that of op.cit.} 
\begin{equation*}
G_{k}((\theta ,\phi ),\tau )=\sum_{\substack{ m,n\in \mathbb{Z}  \\ %
(m,n)\neq (0,0)}}\frac{\theta ^{m}\phi ^{n}}{(m\tau +n)^{k}},
\end{equation*}%
for $\theta ,\phi \in U(1)$ with $G_{k}((1,1),\tau ))=G_{k}(\tau )$ \cite%
{MTZ}. $G_{k}((\theta ,\phi ),\tau )$ is not an eigenfunction of $T(n)$ in
general. However, for prime $p$ 
\begin{equation}
T(p)G_{k}((\theta ,\phi ),\tau )=p^{k-1}G_{k}((\theta ,\phi ),\tau
)+R(p)G_{k}((\theta ,\phi ),\tau ).  \label{Tp_Eisen_twist}
\end{equation}%
Hence if $(\theta ,\phi )=(\theta ^{p},\phi ^{p})$ then $T(p)G_{k}((\theta
,\phi ),\tau )=\sigma _{k-1}(p)G_{k}((\theta ,\phi ),\tau )$.

\section{The Permutation Orbifold of a $C=24$ Holomorphic Vertex Operator
Algebra}

\subsection{The Orbifold of a Holomorphic VOA}

We now consider a Vertex Operator Algebra $V$ (VOA) of central charge $24$
e.g. \cite{FLM}, \cite{Ka}, \cite{MN}. We assume that $V$ is a Holomorphic
VOA (HVOA) so that $V$ is the unique irreducible module for itself with
modular invariant meromorphic partition function \cite{Sch}, \cite{DM} 
\begin{eqnarray*}
Z_{V}(\tau ) &=&Tr_{V}(q^{L(0)-1})=\sum\limits_{k\geq -1}a(k)q^{k} \\
&=&J(\tau )+a(0).
\end{eqnarray*}%
For example, the Moonshine Module $V^{\natural }$ is a HVOA with $%
Z_{V^{\natural }}(\tau )=J(\tau )$ whereas $Z_{V_{L}}(\tau )=J(\tau )+24$
for the Leech lattice HVOA $V_{L}$ \cite{FLM}.

Let $G$ be a finite subgroup of the automorphism group of $V$. Then for $%
g\in G$ define the orbifold trace function 
\begin{equation*}
Z_{V}((g,1),\tau )=Tr_{V}(gq^{L(0)-1}),
\end{equation*}%
(so that $Z_{V}((1,1),\tau )=Z_{V}(\tau )$). For the Moonshine Module $%
V^{\natural }$ 
\begin{equation}
T_{g}(\tau )=Z_{V^{\natural }}((g,1),\tau ),  \label{Tg_def}
\end{equation}%
is the McKay-Thompson series for $g\in \mathbb{M}$, the Monster group of
automorphisms of $V^{\natural }$.

Since $V$ is holomorphic, there is a unique twisted module $M_{h}$ for each $%
h\in G$ $\ $\cite{DLM}. For $g\in C(h)$, the $h$ centralizer, $g$ induces a
class of linear maps $\phi (g)$ on $M_{h}$ so that we may define a twisted
orbifold trace function\footnote{%
denoted by $Z(h,g^{-1},\tau )$ in ref. \cite{DLM},} 
\begin{equation}
Z((g,h),q)=Z((g,h),\tau )=Tr_{M_{h}}(\phi (g)q^{L(0)-1}),  \label{Zgh}
\end{equation}%
a meromorphic function for $\tau \in \mathbb{H}$ [op.cit.]. We define a
right action of the modular group for $\gamma \in \Gamma $ as follows%
\footnote{%
The $0$ subscript denotes the modular weight of $Z((g,h),\tau )$.} 
\begin{equation}
(Z|_{0}\gamma )((g,h),\tau )=Z(\gamma (g,h),\gamma \tau ),  \label{Zghgamma}
\end{equation}%
with 
\begin{equation*}
\gamma (g,h)=(g^{a}h^{b},g^{c}h^{d}).
\end{equation*}%
The trace function enjoys the modular invariance property \ 

\begin{equation}
(Z|_{0}\gamma )((g,h),\tau )=\epsilon _{\gamma }(g,h)Z((g,h),\tau ).
\label{Zghmod}
\end{equation}%
for cocycle $\epsilon _{\gamma }(g,h)\in \mathbb{C}^{\ast }$ [op.cit.]
generalizing earlier ideas of Zhu concerning trace functions \cite{Z}.
Specializing to the McKay-Thompson series (\ref{Tg_def}), these results
imply that $T_{g}(\tau )$ is a meromorphic function on $\mathbb{H}$
satisfying the modular invariance property (\ref{Zghmod}) for $h=1$.

Let us consider orbifolds without a global phase anomaly i.e. where each $%
\phi (g)$ acting on $M_{h}$ can be chosen such that $\epsilon _{\gamma
}(g,h)=1$ so that \cite{Va} 
\begin{equation}
(Z|_{0}\gamma )((g,h),\tau )=Z((g,h),\tau ).  \label{anomaly}
\end{equation}%
In particular, this condition implies that for $h$ of order $m$%
\begin{equation}
Z((1,h),\tau )=\sum\nolimits_{k\in \mathbb{Z}}a((1,h),\frac{k}{m})q^{k/m},
\label{Z1h}
\end{equation}%
for some integers $a((1,h),\frac{k}{m})\geq 0$. We next define the $G-$%
orbifold partition function by

\begin{eqnarray}
Z^{G\mathrm{-orb}}(\tau ) &=&\frac{1}{\left\vert G\right\vert }\sum\limits 
_{\substack{ g,h\in G  \\ gh=hg}}Z((g,h),\tau )  \notag \\
&=&\sum\limits_{[h]\in G}\frac{1}{\left\vert C(h)\right\vert }%
\sum\limits_{g\in C_{G}(h)}Z((g,h),\tau ),  \label{ZGorb}
\end{eqnarray}%
for centralizer $C(h)=\{hg=gh|g\in G\}$ and where $[h]$ denotes a conjugacy
class of $G$. Clearly $Z^{G\mathrm{-orb}}(\tau )$ is also modular invariant.
The most well-known example is the original construction for the Moonshine
Module $V^{\natural }$ as a $\mathbb{Z}_{2}$ orbifold of the Leech lattice
VOA \cite{FLM}.

\subsection{Permutation Orbifolds}

Let $V^{\otimes n}=V\otimes V\otimes \ldots V$ denote the $n^{th}$ tensor
product VOA with partition function $Z_{V^{\otimes n}}(\tau )=Z(\tau )^{n}$
and central charge $24n$. The symmetric group $S_{n}$ naturally acts on $%
V^{\otimes n}$ as an automorphism group. For each $\beta \in S_{n}$ there is
a unique $\beta $-twisted $V^{\otimes n}$ module $M_{\beta }$ which can be
explicitly constructed from the original HVOA $V$ \cite{DMVV}, \cite{BDM}, 
\cite{Ba}. Furthermore, we may explicitly compute the permutation orbifold
partition function $Z^{S_{n}-\mathrm{orb}}(\tau ).$

We illustrate this in the first non-trivial case for $V\otimes V$. Then $%
S_{2}=\langle \sigma \rangle $ for $2-$cycle $\sigma $ where $\sigma
:u\otimes v\rightarrow v\otimes u$ for all $u\otimes v\in V\otimes V$. We
thus find 
\begin{equation*}
Z_{V\otimes V}((\sigma ,1),\tau )=Z(2\tau ).
\end{equation*}%
The $\sigma -$twisted module $M_{\sigma }$ has partition function

\begin{equation*}
Z_{V\otimes V}((1,\sigma ),\tau )=Tr_{M_{\sigma }}(q^{L(0)-1})=Z(\frac{\tau 
}{2}),
\end{equation*}%
with 
\begin{equation*}
Z_{V\otimes V}((\sigma ,\sigma ),\tau )=Z(\frac{\tau +1}{2}),
\end{equation*}%
following (\ref{anomaly}). Thus we obtain%
\begin{equation}
Z^{S_{2}\mathrm{-orb}}(\tau )=\frac{1}{2}Z(\tau )^{2}+T(2)Z(\tau ),
\label{ZS2orb}
\end{equation}%
for Hecke operator $T(2)$ of (\ref{Tnf}) for weight zero.

In general, for $\beta \in $ $S_{n}$, consider the cycle decomposition%
\begin{equation}
\beta =\sigma _{1}^{m_{1}}\sigma _{2}^{m_{2}}\ldots \sigma _{n}^{m_{n}},
\label{gcycle}
\end{equation}%
where $\sigma _{k}$ denotes a $k-$cycle. The conjugacy classes of $S_{n}$
are enumerated by the set of partitions of $n=\sum\limits_{1\leq k\leq
n}km_{k}$. The centralizer is then%
\begin{equation}
C(\beta )=S_{m_{1}}\times (S_{m_{2}}\rtimes C_{2}{}^{m_{2}})\times \ldots
\times (S_{m_{n}}\rtimes C_{n}{}^{m_{n}}),  \label{g_Cent}
\end{equation}%
of order $\prod\limits_{1\leq k\leq n}k^{m_{k}}m_{k}!$ with cyclic group $%
C_{k}=\langle \sigma _{k}\rangle $ and $S_{m_{k}}$ the permutation group on
the $m_{k}$ cycles $\sigma _{k}$. We may construct $M_{\beta }=\otimes
_{k}M_{\sigma _{k}}^{\otimes m_{k}}$ which has partition function \cite{BDM}

\begin{equation*}
Z((1,\beta ),\tau )=Tr_{M_{\beta }}(q^{L(0)-1})=\prod\limits_{1\leq k\leq
n}Z(\frac{\tau }{k})^{m_{k}}.
\end{equation*}%
One eventually finds that the $S_{n}$ permutation orbifold partition
function is \cite{DMVV} 
\begin{equation}
Z^{S_{n}\mathrm{-orb}}(\tau )=\sum\limits_{[\beta ]\in S_{n}}\frac{1}{%
\left\vert C(\beta )\right\vert }\sum\limits_{\alpha \in C(\beta )}Z((\alpha
,\beta ),\tau )=\sum\limits_{\substack{ m_{1},\ldots m_{n}  \\ \sum km_{k}=n
}}\prod\limits_{1\leq k\leq n}\frac{1}{m_{k}!}(T(k)Z(\tau ))^{m_{k}},
\label{ZSnOrb}
\end{equation}%
for the classical Hecke operator $T(k)$ of (\ref{Tnf}).

It is natural to define a permutation orbifold generating function by 
\begin{equation}
Z^{\mathrm{perm}}(p,q)=1+\sum\limits_{n\geq 1}p^{n}Z^{S_{n}\mathrm{-orb}%
}(\tau ),  \label{Zpermgen}
\end{equation}%
for a formal parameter $p$. Thus we obtain [op.cit.] 
\begin{equation}
Z^{\mathrm{perm}}(p,q)=\exp (\sum\limits_{n\geq 1}p^{n}T(n)Z(\tau
))=\prod\limits_{r\geq 1,s\in \mathbb{Z}}\frac{1}{(1-p^{r}q^{s})^{a(rs)}},
\label{Zpermpq}
\end{equation}%
where $Z(\tau )=\sum\nolimits_{k\geq -1}a(k)q^{k}$. This is clearly of the
form\ of the \textbf{inverse} of the LHS of denominator formula (\ref%
{denom_formula}). Thus such expressions canonically arise in the context of
permutation orbifolds for $C=24$ HVOAs.

Specializing to the case of the Moonshine module $V^{\natural }$ where $%
Z(\tau )=J(\tau )$ we obtain 
\begin{gather*}
Z_{V^{\natural }}^{\mathrm{perm}}(p,q)=\exp (\sum\limits_{n\geq
1}p^{n}T(n)J(\tau ))=\frac{1}{pJ(p)-pJ(q)} \\
=1+pJ\left( q\right) +\left( J\left( q\right) ^{2}-c\left( 1\right) \right) {%
p}^{2}+\left( J\left( q\right) ^{3}-2\,J\left( q\right) c\left( 1\right)
-c\left( 2\right) \right) {p}^{3}+ \\
\left( J\left( q\right) ^{4}-3\,c\left( 1\right) J\left( q\right)
^{2}-2\,J\left( q\right) c\left( 2\right) -c\left( 3\right) +c\left(
1\right) ^{2}\right) {p}^{4}+\ldots
\end{gather*}%
This formula and the infinite product formula of (\ref{Zpermpq}) very
strongly suggest that $Z_{V^{\natural }}^{\mathrm{perm}}(p,q)$ is the
partition function for a doubly graded \emph{symmetric }bosonic module with
Monster characters which is, algebraically speaking, the inverse of the
alternating homological structure constructed by Borcherds \cite{Bo2}.
Furthermore, infinite product formulas such as that of (\ref{Zpermpq}) have
been given the interesting interpretation as a "second quantized" string
partition function in the physics literature \cite{DMVV}. A rigorous VOA
construction for such a\ structure would be of obvious interest.

\section{Finite Group and Permutation Orbifolds}

Let us now consider orbifolding $V^{\otimes n}$ with respect to $G\times
S_{n}$ where $G$ acts diagonally on $V^{\otimes n}$ and $S_{n}$ is the
permutation group for $V^{\otimes n}$. We consider again a $C=24$
holomorphic VOA, with modular invariant partition function and where (\ref%
{anomaly}) is holds. We may construct unique twisted sectors for each $%
(h,\beta )\in G\times S_{n}$ \cite{BDM} to find \cite{T4} 
\begin{equation}
Z^{S_{n}\mathrm{-orb}}((g,h),\tau )=\sum\limits_{\substack{ m_{1},\ldots
m_{n}  \\ \sum km_{k}=n}}\prod\limits_{1\leq k\leq n}\frac{1}{m_{k}!}%
(T(k)Z((g,h),\tau ))^{m_{k}},  \label{Snghorb}
\end{equation}%
where here $T(k)$ is the twisted Hecke Operator\ of (\ref{Tnf_twist}). Then (%
\ref{Tnf_gam_twist}) implies 
\begin{equation}
(T(n)Z((g,h),\tau ))|_{0}\gamma =T(n)Z((g,h),\tau ).  \label{Tn_Gam_Orb}
\end{equation}

We may define a permutation orbifold generating function generalizing (\ref%
{Zpermgen}) as follows: 
\begin{eqnarray}
Z^{\mathrm{perm}}((g,h),p,q) &=&\sum\limits_{n\geq 1}p^{\frac{n}{m}}Z^{S_{n}%
\mathrm{-orb}}((g,h),\tau )  \notag \\
&=&\exp (\sum\limits_{n\geq 1}p^{\frac{n}{m}}T(n)Z((g,h),\tau )).
\label{Zpermgengh}
\end{eqnarray}%
where $h$ is of order $m$ \ and using (\ref{Snghorb}). For $g=1$, this
expression reduces to an infinite product formula generalizing (\ref{Zpermpq}%
) to find 
\begin{equation}
Z^{\mathrm{perm}}((1,h),p,q)=\prod\limits_{r\geq 1,s\in \mathbb{Z}}(1-p^{%
\frac{r}{m}}q^{\frac{s}{m}})^{-a((1,h^{r}),\frac{rs}{m})},
\label{Zperm(1,h)}
\end{equation}%
where $Z((1,h^{r}),\tau ))=\sum\nolimits_{k\in \mathbb{Z}}a((1,h^{r}),\frac{k%
}{m})q^{\frac{k}{m}}$.

\section{Monstrous and Generalized Moonshine - the Genus Zero Property}

We now consider the FLM\ Moonshine Module VOA $V^{\natural }$ \cite{FLM} and
its relationship to Moonshine. The original Monstrous Moonshine paper of
Conway and Norton described evidence for an unexpected relationship between
properties of the Monster finite group and the theory of modular forms \cite%
{CN}. Many of these relationships are now understood to be generic to
orbifold constructions in conformal field theory/VOA theory e.g. \cite{FLM}, 
\cite{T1}, \cite{T2}, \cite{DLM}. However, the special feature that sets the
Moonshine Module $V^{\natural }$ apart from other VOAs is the Genus Zero
Property \cite{CN}. This states that for each $g\in \mathbb{M}$, the
McKay-Thompson series $T_{g}(\tau )$ of (\ref{Tg_def}) is a hauptmodul for
some genus zero modular group $\Gamma _{g}$. Thus for $g$ of prime order $%
o(g)=$ $p$, one finds (excluding one class of order 3) that either $\Gamma
_{g}=\Gamma _{0}(p)$ with $g=p-$ (in the notation of \cite{CN}) or else with 
$g=p+$ with $\Gamma _{g}=\Gamma _{0}(p)+=\langle \Gamma _{0}(p),W_{p}\rangle 
$ where $\Gamma _{0}(p)=\{\left( 
\begin{array}{cc}
a & b \\ 
c & d%
\end{array}%
\right) |c=0\ \mathrm{mod}\ p\}$ and $W_{p}:\tau \rightarrow -1/p\tau $ is a
Fricke involution. In general, we say that $g\in \mathbb{M}$ is Fricke if $%
T_{g}$ is invariant under a Fricke involution $W_{N}:\tau \rightarrow
-1/N\tau $ where $N=ko(g)$ and $k|24$ and is otherwise non-Fricke. $k=1$ in
the global phase anomaly free cases where (\ref{anomaly}) holds.

The distinction between Fricke and non-Fricke classes is particularly
important in the orbifold interpretation of Moonshine \cite{T1}. There is
very significant evidence for the general conjecture that the genus zero
property for a McKay-Thompson series is equivalent to the statement that for
any global phase anomaly free element $g$, orbifolding $V^{\natural }$ with
respect to $\langle g\rangle $ for $g$ Fricke results in $V^{\natural }$
again whereas orbifolding $V^{\natural }$ with respect to $\langle g\rangle $
for $g$ non-Fricke results in the Leech lattice VOA \cite{T2}.

Generalized Moonshine refers to the still generally unproven conjecture of
Norton \cite{N2} that for each commuting pair $g,h\in \mathbb{M}$, then $%
Z_{V^{\natural }}((g,h),\tau )$ is either a hauptmodul for a genus zero
modular group or is a constant. \ It is easy to show using (\ref{Zghmod})
that (1) $Z_{V^{\natural }}((g,h),\tau )$ is constant iff $g^{c}h^{d}$ is
non-Fricke for all $(c,d)=1$ \cite{N2}, \cite{T3} and (2) if $g,h\in \langle
k\rangle $ for some $k\in \mathbb{M}$ then $Z_{V^{\natural }}((g,h),\tau )$
is a hauptmodul (since it can then be modular transformed to a
McKay-Thompson series \cite{T3}, \cite{IT1}, \cite{DLM}). \ In the remaining
"non-trivial" cases, we may use (\ref{Zghmod}) again to transform $%
Z_{V^{\natural }}((g,h),\tau )$ to a trace over a Fricke twisted module so
that the genus zero property reduces to an analysis of $h$ Fricke cases
alone. The case of $h=2+$, with centralizer $2.B$ for the Baby Monster $B$,
has now been proved by Hoehn \cite{H}. The relationship between orbifoldings
of the Moonshine module and the genus zero property of Generalized Moonshine
for $h=p+$ is discussed at length in \cite{T3}, \cite{IT1}, \cite{IT2}, \cite%
{I}.

The approach taken in Borcherds' proof of the Monstrous Moonshine genus zero
property is to firstly prove a twisted denominator identity generalizing (%
\ref{denom_formula}) \cite{N1}, \cite{Bo2}%
\begin{equation}
\exp (-\sum\limits_{n\geq 1}p^{n}T(n)T_{g}(\tau ))=p(T_{g}(p)-T_{g}(q)).
\label{denom_Tg}
\end{equation}%
This is the defining formula for completely replicable functions \cite{N1}, 
\cite{FMN} from which it follows that the leading coefficients of $%
T_{g^{i}}(\tau )$ for $i=1,2,\ldots $ determine $T_{g}(\tau )$. Koike showed
that the list of hauptmoduln appearing in the Moonshine Conjectures are
completely replicable \cite{Ko}. Based on this result and an analysis of the
leading coefficients (using the explicit form for $T_{g}(\tau )$ found by
FLM for $2-$ centralizers in the Monster \cite{FLM}) Borcherds then
demonstrated that indeed $T_{g}(\tau )$ obeys the  genus zero property. This
part of the proof was improved upon in \cite{CG} where meromorphicity,
modularity and the genus zero property are shown to generally follow from
the infinitely many replication formulas that follow from (\ref{Pnf}) and (%
\ref{denom_Tg}), namely%
\begin{equation}
T(n)T_{g}(q)=\frac{1}{n}F_{n}(T_{g}(q)),  \label{Replic_Tg}
\end{equation}%
where $F_{n}$ is the Faber polynomial for $T_{g}(q)$.

However, as already noted, $T_{g}(q)$ is known to be meromorphic on $\mathbb{%
H}$ from \cite{DLM}. Furthermore, by combining (\ref{Replic_Tg}) with the
general modular transformation property (\ref{Zghmod}) (or (\ref{anomaly})
in the absence of global anomalies) one can also expect to obtain the genus
zero property for $T_{g}$ in a more direct fashion along the lines of the
methods described in refs. \cite{T1}, \cite{T2}.

On the other hand, in the permutation orbifold construction based on $%
V^{\natural }$, (\ref{denom_Tg}) reads%
\begin{eqnarray*}
Z_{V^{\natural }}^{\mathrm{perm}}((g,1),p,q) &=&\exp (\sum\limits_{n\geq
1}p^{n}T(n)Z_{V^{\natural }}((g,1),q)) \\
&=&\frac{1}{p(Z_{V^{\natural }}((g,1),p)-Z_{V^{\natural }}((g,1),q))}.
\end{eqnarray*}%
Recall our previous remarks concerning the reduction of Generalized
Moonshine to Fricke classes. Consider $h$ a global phase anomaly-free Fricke
element of order $o(h)=m$ so that $Z_{V^{\natural
}}((1,h),q^{m})=Z_{V^{\natural }}((h,1),q)$. It follows that%
\begin{equation*}
\exp (\sum\limits_{n\geq 1}p^{n}T(n)Z((1,h),q^{m}))=\frac{1}{%
p(Z_{V^{\natural }}((1,h),p^{m})-Z_{V^{\natural }}((1,h),q^{m}))}.
\end{equation*}%
Hence from (\ref{Zperm(1,h)}) we find for such Fricke $h$ that

\begin{equation}
Z_{V^{\natural }}^{\mathrm{perm}}((1,h),p,q)=\frac{1}{p^{\frac{1}{m}%
}(Z_{V^{\natural }}((1,h),p)-Z_{V^{\natural }}((1,h),q))}.
\label{Zperm(1,h)moonshine}
\end{equation}%
This together with the general result (\ref{Zperm(1,h)}) again suggests the
existence of a symmetric bosonic construction forming a doubly-graded module
for the centralizer of $C(h)$ for each such Fricke element $h$. It is thus
natural to conjecture that for all order $m$ global phase anomaly-free
Fricke elements $h$ the following holds 
\begin{eqnarray}
Z_{V^{\natural }}^{\mathrm{perm}}((g,h),p,q) &=&\exp (\sum\limits_{n\geq
1}p^{\frac{n}{m}}T(n)Z_{V^{\natural }}((g,h),q))  \notag \\
&=&\frac{1}{p^{\frac{1}{m}}(Z_{V^{\natural }}((g,h),p)-Z_{V^{\natural
}}((g,h),q))}.  \label{GenReplication}
\end{eqnarray}%
From (\ref{Pnf}), this is equivalent to the following Generalized Moonshine
replication formula for global phase anomaly-free Fricke elements $h$ 
\begin{equation}
T(n)Z_{V^{\natural }}((g,h),q)=\frac{1}{n}F_{n}(Z_{V^{\natural }}((g,h),q)),
\label{Replic_Tgh}
\end{equation}%
where $F_{n}$ is the Faber polynomial for $Z_{V^{\natural }}((g,h),q^{m}))$.
The equivalence of this replication formula to the genus zero property for
Generalized Moonshine would therefore require a suitable generalization of
the various Monstrous Moonshine arguments.

We conclude with the example of $n=2$. Then (\ref{Replic_Tgh}) implies 
\begin{eqnarray*}
&&Z_{V^{\natural }}((g^{2},h),2\tau )+Z_{V^{\natural }}((g,h^{2}),\frac{\tau 
}{2})+Z_{V^{\natural }}((gh,h^{2}),\frac{\tau +1}{2}) \\
&=&Z_{V^{\natural }}((g,h),\tau )^{2}-2a((g,h),\frac{1}{m}),
\end{eqnarray*}%
where $Z_{V^{\natural }}((g,h),\tau )=$ $q^{-\frac{1}{m}}+0+a((g,h),\frac{1}{%
m})$ $q^{\frac{1}{m}}+...$ which corresponds to an example quoted in ref. 
\cite{N2}. In particular, for $h=2+$ and $g$ of order $2\ $ then using
Fricke invariance we find%
\begin{equation*}
T_{h}(\tau )+T_{g}(\frac{\tau }{2})+T_{gh}(\frac{\tau +1}{2})=Z_{V^{\natural
}}((g,h),\tau )^{2}-2a((g,h),\frac{1}{2}),
\end{equation*}%
which can be easily verified in each case.

\end{document}